\documentstyle[psfig]{article}

\input mssymb12.tex
\def\NABLA#1{{\mathop{\nabla\kern-.5ex\lower1ex\hbox{$#1$}}}}
\def\Nabla#1{\nabla\kern-.5ex{}_#1}

\newcommand{\p}{{\partial}}

\newtheorem{theorem}{Theorem} 
\newtheorem{corollary}[theorem]{Corollary}
\newtheorem{lemma}[theorem]{Lemma}

\title{Lagrangian isotopies in Stein manifolds}
\author{R. Hind\thanks{Supported in part by NSF grant DMS-0204634.}}

\date{\today}

\begin{document}

\maketitle

\section{Introduction}

Studying the space of Lagrangian submanifolds is a fundamental
problem in symplectic topology. Lagrangian spheres appear
naturally in the Leftschetz pencil picture of symplectic
manifolds.

In this paper we demonstrate the uniqueness up to Hamiltonian
isotopy of the Lagrangian spheres in some $4$-dimensional Stein
symplectic manifolds. The most important example is the cotangent
bundle of the $2$-sphere, $T^* S^2$, with its standard symplectic
structure.

We recall that if a convex symplectic manifold has a boundary
of contact-type, then we can perform surgery operations on
the manifold by adding handles to the boundary. In the $4$-dimensional
case these handles can be of index $1$ or $2$. Our other examples
are symplectic manifolds formed by adding $1$-handles to a
unit cotangent bundle $T^1 S^2$. Questions regarding Lagrangian
isotopy classes are independent of which metric we use to
define a unit tangent bundle or of any choices involved in
adding $1$-handles.

\begin{theorem}
Let $M$ be $T^* S^2$ or the result of adding any number of
$1$-handles to $T^1 S^2$ and $L \subset M$ be a Lagrangian sphere.
Then there exists a Hamiltonian diffeomorphism of $M$ mapping $L$
onto the zero-section.
\end{theorem}

The statement is false even if we add to $T^1 S^2$ a single
$2$-handle along the Legendrian curve in a single fiber of the
boundary. The result of this surgery operation is a plumbing $M$
of two copies of $T^1 S^2$. The symplectic manifold $M$ has two
Lagrangian spheres $L_1$ and $L_2$ coming from the zero-sections
in the $T^1 S^2$. Now, associated to any Lagrangian sphere $L$ is
a compactly supported symplectomorphism $\tau_L$ (a generalized
Dehn Twist) whose square is smoothly but not necessarily
symplectically isotopic to the identity. Thus $\tau^2_{L_2}(L_1)$
is a Lagrangian sphere in $M$ which is smoothly isotopic to $L_1$.
However, as demonstrated by P. Seidel in \cite{seidel}, a Floer
homology computation shows that $\tau^2_{L_2}(L_1)$ is not
Hamiltonian isotopic to $L_1$.

We will establish the theorem by utilizing an existence result
for almost-complex structures on $S^2 \times S^2$ with convenient
properties, taken from \cite{hind}, and a fact about diffeomorphisms
of the $2$-sphere.

In section $2$ we prove our result on the diffeomorphisms of $S^2$.
In section $3$, by using the conclusions of \cite{hind}, we reduce
our theorem in the case of $M=T^* S^2$ to the statements in section
$2$. In section $4$ we deal with the addition of handles. This
involves slightly generalizing the results from \cite{hind} so
we will review them again there.

As yet we are unable to prove any similar results for higher genus
Lagrangian surfaces, but we can make some remarks about the case
of $\Bbb R P^2$. The results of \cite{hind} show that any
Lagrangian sphere $L$ in $S^2 \times S^2$ homotopic to the
antidiagonal $\overline{\Delta}$ is in fact Lagrangian isotopic to
$\overline{\Delta}$. In this paper we will show that if $L$ is
disjoint from the diagonal $\Delta$ then the Lagrangian isotopy
can be chosen to lie in $S^2 \times S^2 \setminus \Delta$. Now,
the involution $\sigma$ of $S^2 \times S^2$ interchanging the two
factors has fixed-point set equal to $\Delta$ and restricts to the
antipodal map on $\overline{\Delta}$. If $L$ is invariant under
$\sigma$ then the isotopy can also be chosen to be
$\sigma$-equivariant. Now, quotienting out by $\sigma$, we observe
that $S^2 \times S^2 \setminus \Delta$ is a double-cover of a unit
cotangent bundle of $\Bbb R P^2$ and Lagrangian spheres in $T^*
\Bbb R P^2$ homotopic to the zero-section therefore correspond to
$\sigma$-invariant Lagrangian spheres in $S^2 \times S^2 \setminus
\Delta$ homotopic to $\overline{\Delta}$. Hence we have the
following corollary.

\begin{corollary}

A Lagrangian $\Bbb R P^2$ homotopic to the zero-section in $T^*
\Bbb R P^2$ must be Hamiltonian isotopic to the zero-section.
\end{corollary}

A natural compactification of the (unit) cotangent bundle of $\Bbb
R P^2$ is $\Bbb C P^2$. The uniqueness up to Hamiltonian isotopy
of Lagrangian $\Bbb R P^2$s inside $\Bbb C P^2$ can be established
by other methods. For example, the surgery technique described in
\cite{sym} replaces a Lagrangian $\Bbb R P^2$ by a symplectic
sphere, transforming $\Bbb C P^2$ into $S^2 \times S^2$. But the
symplectic spheres in $S^2 \times S^2$ have been classified up to
Hamiltonian isotopy in \cite{tian}.

\section{Diffeomorphisms of the two-sphere}

In this section we let $f$ denote a diffeomorphism of the $2$-sphere
$S^2$ and for a point $x\in S^2$ we denote its antipodal point by
$-x$.

We say that a diffeomorphism $f$ has the property $(*)$ if
$f(x) \neq -x$ for all $x \in S^2$.

The aim of the section is to prove the following theorem.

\begin{theorem}
Suppose that $f$ has the property $(*)$. Then there exists an
isotopy $f_t$, $0\le t \le 1$, with $f_0 = \mathrm{id}$ and $f_1=f$
such that $f_t$ has property $(*)$ for all $t$.
\end{theorem}

{\bf Proof of theorem}

Let $E$ denote an equator on $S^2$. The complement of $E$ consists
of two open disks $H_1$ and $H_2$ with $-H_1 =H_2$.

We observe that any diffeomorphism $g$ with property $(*)$ and
which preserves $E$ is indeed isotopic to the identity through
diffeomorphisms $G_t$ also satisfying $(*)$. To construct such an
isotopy, we first isotope $g$ to the identity in a neighbourhood
of $E$. Now the resulting map restricts to a compactly supported
diffeomorphism of $H_1$ and $H_2$. But compactly supported
diffeomorphisms of the disk are isotopic to the identity (see for
instance \cite{thurston}, page 205). Combining these isotopies we
get the required isotopy of $g$. It satisfies $(*)$ since $-H_1
=H_2$.

Hence it suffices to find a suitable isotopy from $f$ to a
diffeomorphism preserving an equator $E$.

Now, the diffeomorphism is necessarily orientation preserving and
so has a fixed point $p_0$. After a small perturbation we may
assume that $f=\mathrm{id}$ near $p_0$. We first assume that $f$
also satisfies $f(-p_0)=-p_0$, and so by another small
perturbation can suppose that $f(x)=x$ also near $-p_0$. Then we
can identify $S^2 \setminus \{p_0,-p_0\}$ with $\Bbb R \times S^1$
such that in these coordinates the antipodal map is given by
$$-(x,e^{i\theta})=(-x,-e^{i\theta})$$
and $f$ is equal to the identity on the complement of a compact
set $C$. In these coordinates let the equator $E=\{x=0\}$.

We construct our isotopy by applying the following lemma.

\begin{lemma}
Let $L_s$, $-\infty <s < \infty$ be a $1$-dimensional foliation of
$\Bbb R \times S^1$ coinciding with $\{ x=s \}$ outside $C$. For
any $N$ there exists a compactly supported isotopy $f_t$
satisfying $(*)$ such that $f_0 =f$ and $f_t(z) \in f(L_{s+tN})$
for all $z \in L_s \cap C$, $0 \le t \le 1$.
\end{lemma}

Given this, we can conclude as follows. First we apply the lemma
with $L_s = \{x=s\}$ to find an isotopy to a new diffeomorphism,
still denoted by $f$, such that $f(z)\in L_K$ for all $z \in E$,
for $K$ large. Let $N=f^{-1}(E)$. Then $N$ is disjoint from $E$
and so we can form another foliation $L_s'$ which includes the
circles $N$ and $E$. After another application of the lemma,
perhaps now with a larger value of $K$, we can isotope $f$ to
another diffeomorphism now satisfying $f(z) \in E$ for all $z \in
E$.

{\bf Proof of lemma}

As the condition on our isotopy is an open one, we may assume any
necessary genericity properties for the diffeomorphism $f$ with respect
to the foliation $L_s$.

Suppose that $N>0$. We can set $f_t
(x,e^{i\theta})=(x+Nt,e^{i\theta})$ for $x$ large and positive. In
general, for $r \in \Bbb R$, let $a_r$ be a diffeomorphism of
$\Bbb R \times S^1$ such that $a_r(L_s)=L_{s+r}$. Outside of $C$
we may assume that $a_r(x,e^{i\theta})=(x+r,e^{i\theta})$. Then we
will define $f_t$ by
$$f_t(z)=h_tfa_{Nt}(z)$$
where $h_t$ is a diffeomorphism of $\Bbb R \times S^1$ which
preserves the foliation $\{L_s \}$. We set
$h_{t,s}=h_t|_{f(L_{s+Nt})}$.

Then we need to find smoothly varying $h_{t,s}$ such that
$h_{t,s}(f(a_{tN}(z)) \neq -z$ for all $z\in L_s$, $s$ and $0 \le t \le 1$.

For $s$ very large we have $h_{t,s} =\mathrm{id}$ and it is
required to show that we can extend these diffeomorphisms for all
parameters $s$.
Again since property $(*)$ is an open condition, we observe that
once we have defined the $h_{t,s_0}$ for some $s_0$ we can smoothly
extend the functions to define $h_{t,s}$ for $s$ slightly less than
$s_0$.

Another observation is that an isotopy defined with the required
properties for $s \ge s_1$, some $s_1$, can always be extended to
an isotopy of $\Bbb R \times S^1$ satisfying the property $(*)$
and preserving the foliation by levels $f(L_s)$. To do this, we
simply smoothly cut-off the vector field generating $f_t$ in such
a way that the cut-off still preserves the levels. Property $(*)$
still holds provided that the vector field is zero on $f(L_s)$ for
$s<s_1 -\epsilon$, $\epsilon$ sufficiently small. The key point is
that since the orbit of a point $x$ on $L_{s_1}$ avoids the point
$-x$, it also avoids $-y$ for all $y$ close to $x$.

For any $s$, as $t$ increases from $0$ to $1$ there is a varying
collection of points $I_{t,s}=f(L_{s+tN})\cap -L_s$. The diffeomorphisms
$h_{t,s}$ can be extended arbitrarily once they are defined on
these intersections. We notice that $I_{t,s}$ is empty for $s$
sufficiently large or small.

Now, for typical $s'$ the families of points $I_{t,s}$, $0 \le t \le 1$,
will vary continuously with $s$ for $s$ close to $s'$. For a fixed value
of $s$ they will consist of a continuously varying set of points which
at certain times appear or vanish in pairs. Continuous variation with
$s$ means that we have smooth families of diffeomorphisms
$$\phi_s:[0,1] \to [0,1]$$
$$g_{t,s}:f(L_{s'+tN}) \to f(L_{s+\phi_s(t)N})$$
$$b_{t,s}:L_{s'} \to L_{s}$$
such that $\phi_{s'}=\mathrm{id}$, $g_{t,s'}=\mathrm{id}$,
$b_{t,s'}=\mathrm{id}$ and $g_{0,s}(f(z))=f(b_{0,s}(z))$.
They can be chosen such that
$g_{t,s}(I_{t,s'})=I_{\phi_s(t),s}$ and if $-z \in I_{t,s'}$,
then $g_{t,s}(-z)=-b_{t,s}(z)$. The existence of such diffeomorphisms
implies that if we have defined suitable $h_{t,s}$ for an $s$ close
to $s'$ then we may define the $h_{t,s'}$ by
$$h_{t,s'}(f(a_{tN}(z))) = g^{-1}_{t,s}(h_{\phi_s(t),s}(f(a_{\phi_s(t)N}b_{t,s}(z))).$$

At a finite collection of parameters $s_i$ the pattern of intersections
$I_{t,s_i}$ will change from nearby values. Assuming $f$ to be generic,
the change will occur only near a single point in $C$ for a single
$t$ parameter.

Suppose that $s'$ is such a critical parameter. In the first case we
consider the situation when $f(L_{s'})$ and $f(L_{s'+N})$ are transverse
to $-L_{s'}$. Then since $s'$ is critical there exists a $\sigma$ with
$s'< \sigma <s'+N$ with $f(L_{\sigma})$ tangent to order three with
$-L_{s'}$. For, if these tangencies are all of order at most two, then
the intersections $I_{t,s}$ will indeed vary continuously with $s$ for
$s$ close to $s'$. Genericity allows us to assume that such a tangency
has order no more than three.

Let the tangency occur at a point $p$. We can choose coordinates $(x,y)$
in $\Bbb R^2$ with $p=(0,0)$ such that the foliation $-L_{s}$ is
given by $\{x=s'-s\}$ and the foliation $\{f(L_s)\}$ by
$\{x=y^3-sy \pm s\}$, where here we take the parameter $\sigma =0$.

If the second foliation is $x=y^3 -sy +s$ then $f(L_s) \cap -L_{s_1}$
is a single point for all $s$ close to $0$ and $s_1 \ge s'$ but
$f(L_s) \cap -L_{s_2}$ consists of three points for some $s>0$ and
$s_2<s'$. In this case we can define the diffeomorphisms $g_{t,s}$
and $b_{t,s}$ as before for $s>s'$ and thus define $h_{t,s'}$. Recall that
we are implicitly assuming that the $h_{t,s}$ can de defined for
$s >s'$ but not necessarily for $s \le s'$.

However, if the second foliation is given by $y=x^3 -sx -s$,
then $f(L_s) \cap -L_{s_1}$ consists of three points for some
$s>0$ and $s_1 >s'$ but only one point for all $s$ close to $0$
and $s_1 \le s'$. Thus the diffeomorphisms $g_{t,s}$ do not exist
as before for $s>s'$. Nevertheless we can define diffeomorphisms
$g_{t,s}$ satisfying the same conditions as before away
from a neighbourhood of $t=\sigma - s'$ and the point $p$.
Near $p$ and $t=\sigma -s'$ we extend the diffeomorphisms sending
$I_{t,s'}$ to the first point in $I_{\phi_s(t),s}$, coherently
ordering the intersection points along $-L_s$. This works as before
to extend the $h_{t,s}$.

Finally suppose that $-L_{s'}$ is tangent to $f(L_{s'})$ or
$f(L_{s'+N})$. We may now assume that this tangency is of
second order and the diffeomorphisms $g_{t,s}$ and $b_{t,s}$
will exist as before for $t$ away from $0$ and $1$. This is
already enough in the case of a tangency with $f(L_{s'})$
since we can set $h_{t,s'}=\mathrm{id}$ for $t$ close to $0$.

If $f(L_{s'+N})$ has a second order tancency with $-L_{s'}$ at
a point $p$ then $f(L_{s+N}) \cap -L_s$ is either empty or
consists of two points for $s>s'$. In the first case it is
easy to extend the $h_{t,s}$ to $s=s'$ simply ensuring that
$f^{-1}_1(-p) \neq p$.

In the second case, suppose that $f_1(L_{s_1}) \cap -L_{s_1} =\{f_1(y),f_1(z)\}$
for some $s_1 >s'$. We can extend the isotopy to $L_{s'}$ if and only
if the points $f(y),-y,-z,f(z)$ do not occur in that order along
$-L_{s_1}$, this is the obstruction to the intersection points moving
together. But recall that the isotopy does indeed extend
with property $(*)$ to all
$L_s$ if we neglect the condition that $f_1(L_s)=f(L_{s+N})$, but still
can require that $f_1(L_s)=f_1(L_r)$ for some $r>s$. The extended
isotopy then still has the property that some $f_1(L_s)$ is tangent
to $-L_s$ near the point $p$. Thus the obstruction must vanish and
we can define $f_t(L_{s'})$ as required. This completes the proof
of the lemma.

\vspace{0.1in}

To complete the proof of Theorem $3$ we must justify the
assumption that $f(-p_0)=-p_0$. Specifically, assuming that
$f=\mathrm{id}$ in a neighborhood $U$ of $p_0$ we need to find an
isotopy of $f$ to a diffeomorphism which is also equal to the
identity near $-p_0$.

Now let $F$ be the equator passing through $p_0$, $-p_0$ and
$f(-p_0) \neq -p_0$, and $G$ be a small circle in $U$ touching $F$
only in a small interval about the point $p_0$. We can extend the
circles $G$ and $F$ to a $1$-parameter family of circles $C_s$,
$-\infty < s < \infty$, with $F=C_0$, $G=C_1$ and $C_s \subset U$
whenever $|s|>1$. Each of the circles will intersect $F$ in a
small interval around $p_0$ and away from $p_0$ they give a
foliation of $S^2$. Then by a variation of Lemma $4$ we may find
an isotopy of $f$ to another diffeomorphism $f$ satisfying
$f(F)=G$ with $f=\mathrm{id}$ still on a small neighborhood of
$p_0$. Consider a path $\gamma$ from $f(-p_0)$ to $-p_0$ following
part $\gamma_1$ of $G$ and then part $\gamma_2$ of the equator $F$
avoiding $p_0$. In fact we may suppose that
 $\gamma_1 =f(\gamma_2)$. We observe that on the path $\gamma$ we will never
 encounter points $f(z)$ then $-z$ in that order. This is clear
 because $\gamma_1=f(\gamma_2)$ and $\gamma_2$ is disjoint from
 $-\gamma_2$. Hence, by following $\gamma$, we can
 find an isotopy $g_t$, supported in a small neighbourhood of $\gamma$ and such that
 $g_0 =\mathrm{id}$ and $g_1(f(-p_0))=-p_0$. Furthermore $g_t
 \circ f$ satisfies property $(*)$ and so gives our isotopy as
 required.

\section{Lagrangian spheres in $T^* S^2$}

Let $L$ be a Lagrangian sphere in $T^* S^2$. This has
self-intersection number $-2$ and so must be homotopic to the
zero-section. By scaling in the fibers we may assume that $L
\subset T^1 S^2$. We will identify $T^1 S^2$ with the complement
of the diagonal $\Delta$ in $S^2 \times S^2$ with its standard
split symplectic form $\omega = \omega_0 \oplus \omega_0$. Under
this identification, the zero-section in $T^1 S^2$ becomes the
antidiagonal $\overline{\Delta}$. Thus our theorem in this case
reduces to the following.

\begin{theorem} Given a Lagrangian sphere $L \subset S^2 \times S^2 \setminus \Delta$ homotopic to $\overline{\Delta}$, there exists a Hamiltonian
isotopy of $S^2 \times S^2$ which fixes $\Delta$ and maps $L$ onto
$\overline{\Delta}$.
\end{theorem}

Given an almost-complex structure $J$ on $S^2 \times S^2$ tamed by $\omega$,
Gromov showed in \cite{gr} that
there exist unique foliations ${\cal F}_0$ and ${\cal F}_1$ by $J$-holomorphic
curves in the classes $[S^2 \times \mathrm{pt}]$ and
$[\mathrm{pt} \times S^2]$. With respect to the standard almost-complex
structure $J_0= i \oplus i$, these foliations are exactly
$S^2 \times \mathrm{pt}$ and $\mathrm{pt} \times S^2$. The key
lemma which we need from \cite{hind} is the following.

\begin{lemma}
There exists a tame almost-complex structure $J$ on $S^2 \times S^2$
such that each curve in the corresponding foliations ${\cal F}_0$
and ${\cal F}_1$ intersects $L$ transversally in a single point.
The almost-complex structure $J$ can be taken to agree with
$J_0$ near $\Delta$.
\end{lemma}

The second statement was not included in \cite{hind} but is clearly
true from the proof.

There exists a family of tame almost-complex structures $J_t$,
$0 \le t \le 1$ on $S^2 \times S^2$ with $J_1 =J$ and, for all
$t$, $J_t =J_0 =i\oplus i$ near $\Delta$. In particular, $\Delta$
is a $J_t$-holomorphic curve for all $t$. By the positivity of
intersections for $J_t$-holomorphic curves, each holomorphic curve
in the foliations ${\cal F}_0$ and ${\cal F}_1$ intersects
$\Delta$ transversally in a single point.

We define a diffeomorphism $f:\Delta \to \Delta$ by $f(x)=y$, where
$y \in \Delta$ is the unique point such that the $J$-holomorphic
curve in ${\cal F}_1$ through $y$ intersects the $J$-holomorphic
curve in ${\cal F}_0$ through $x$ on $L$. Then $f(x) \neq x$ for
all $x \in \Delta$.

As in the previous section, for a point $x \in \Delta$ we denote its
image under the antipodal map by $-x$. Then the $J_0$-holomorphic
curve in ${\cal F}_0$ through $x$ intersects the $J_0$-holomorphic
curve in ${\cal F}_1$ through $-x$ on $\overline{\Delta}$ for all
$x \in \Delta$.

We can apply the theorem of section $2$ to get the following.

\begin{lemma}
There exists an isotopy $g_t:\Delta \to \Delta$, $0 \le t \le 1$,
with $g_0 = \mathrm{id}$, $g_1 = -f^{-1}$ and $g_t(x)\neq -x$ for all
$t$ and $x \in \Delta$.
\end{lemma}

We now define maps $\phi_t :S^2\times S^2 \to S^2 \times S^2$ by
requiring that $\phi_t$ maps the $J_t$-holomorphic curves in
${\cal F}_0$ and ${\cal F}_1$ to the corresponding $J_0$-holomorphic
foliations, the $J_t$-holomorphic curve in ${\cal F}_0$ through $x\in \Delta$ maps
to the $J_0$-holomorphic curve in ${\cal F}_0$ through $x$ and the $J_t$-holomorphic
curve in ${\cal F}_1$ through $x$ maps to the $J_0$-holomorphic
curve in ${\cal F}_1$ through $g_t(x)$.

Then $\phi_0 =\mathrm{id}$, $\phi_1(L)=\overline{\Delta}$ and
$\phi_t(\Delta)$ is disjoint from $\overline{\Delta}$ for all
$t$. Let $L_t= \phi_t^{-1}(\overline{\Delta})$, so $L_t$ gives
a smooth isotopy from $L$ to $\overline{\Delta}$ in
$S^2 \times S^2 \setminus \Delta$.

Also, ${\phi_t}_*(J_t)$ is tamed by the split form $\omega$, and
we see from this that $\phi_t(\Delta)$ is a symplectic submanifold
for all $t$.

For fixed $t$, set $\omega_s=s\phi_t^*(\omega)+(1-s)\omega$. This
is a symplectic form for all $0 \le s \le 1$. It is clearly closed
and is symplectic since it tames $J_t$. We note that $\Delta$ is
symplectic for all $\omega_s$ and, if $t=0$ or $t=1$, $L_t$ is
Lagrangian with respect to all $\omega_s$. Hence by an application
of Moser's theorem we can find a diffeomorphism $\psi_t$ of $S^2
\times S^2$ such that $\psi_t^*(\omega)=\phi_t^*(\omega)$. The
$\psi_t$ can be chosen to vary smoothly with $t$, to fix $\Delta$
and such that $\psi_0=\mathrm{id}$ and $\psi_1$ fixes $L$. To see
this, we recall that Moser's method involves writing
$\omega_s=\omega_0+d\alpha_s$ and studying the flow of the
vectorfield $X_s$ defined by $X_s \rfloor \omega_s
=\frac{d\alpha_s}{ds}$. The definition implies that ${\cal
L}_{X_s} \omega_s =d( \frac{d\alpha_s}{ds})=\frac{d\omega_s}{ds}$.
We have the freedom in this construction to add any smooth family
of exact $1$-forms $\beta_s$ to the $\alpha_s$. These $\beta_s$
can be chosen such that $\alpha_s +\beta_s$ vanishes on the
symplectic normal bundle to $\Delta$ and, if $t=0$ or $t=1$, on
the tangent bundle to $L_t$. Then the flow fixes $\Delta$ and, if
$t=0$ or $t=1$, also fixes $L_t$.

Thus $\psi_t(L_t)$ is a Lagrangian isotopy
from $L$ to $\overline{\Delta}$ inside $S^2 \times S^2 \setminus \Delta$
as required.

\section{Manifolds with $1$-handles}

We now consider the class of convex symplectic manifolds
constructed by adding $1$-handles to the unit cotangent bundle
$T^1 S^2$. Our first observation is that any such manifold $M$ can
be symplectically embedded in $(S^2 \times S^2, \omega)$, after
perhaps scaling the symplectic form. This follows from the methods
of \cite{elg}. We can arrange that the zero-section in $T^1 S^2$
again becomes identified with $\overline{\Delta}$ and the boundary
of $M$ is a smooth hypersurface $\Sigma$ of contact-type in $S^2
\times S^2$.

We plan to find families of almost-complex structures $J_t$ on
$S^2 \times S^2$ and diffeomorphisms $f_t:\Delta \to \Delta$ such
that the $J_t$-holomorphic curves in ${\cal F}_0$ through points
$x$ and the $J_t$-holomorphic curves in ${\cal F}_1$ through
$f_t(x)$ intersect on embedded spheres $L_t \subset M$ with $L_0
=\overline{\Delta}$ and $L_1=L$. The almost-complex structures can
be constructed by deforming $J_0$ in a neighbourhood of $\Sigma$
and, for $t$ close to $0$ or $1$, also in a neighbourhood of
$\overline{\Delta}$ or $L$.

Suppose that we perform the operation of stretching-the-neck along
$\Sigma$. That is, we symplectically identify a neighbourhood of
$\Sigma$ in $S^2 \times S^2$ with $((-\epsilon,\epsilon)\times
\Sigma,d(e^t\alpha))$, where $\alpha$ is a fixed contact form on
$\Sigma$. We can then produce a manifold $A_N$ by replacing this
neighbourhood by $(-N,N)\times \Sigma$. Our original
almost-complex structure can be extended over $(-N,N)\times
\Sigma$ to be translation invariant and the symplectic form can be
extended over $(-N,N)\times \Sigma$ such that $A_N$ is
symplectomorphic to $(S^2 \times S^2, \omega)$ via a
symplectomorphism equal to the identity outside $(-N,N)\times
\Sigma$. Under this symplectomorphism we can think of stretching
the neck as studying a family of almost-complex structures $J_N$
on $S^2 \times S^2$ which degenerate along $\Sigma$ as $N \to
\infty$.

At the same time, we can deform the almost-complex structure along
the boundary of tubular neighborhoods $U_0$ or $U_1$ of $L_0
=\overline{\Delta}$ or $L_1=L$ respectively. Stretching to length
$N$ on both contact hypersurfaces $\Sigma$ and $\p U_i$ we obtain
almost-complex structures $J_{N,0}$ and $J_{N,1}$. There exist
smooth families of almost-complex structures $J_{N,t}$ connecting
$J_{N,0}$ and $J_{N,1}$ which are fixed on the tubular
neighborhoods of $\Sigma$ and in the complement of $M$.

Following the work of Hofer, see \cite{hwz}, and as in
\cite{hind}, after taking suitable subsequences, for $i,j=0,1$
families of $J_{N,i}$-holomorphic curves in ${\cal F}_j$ will
converge to unions of finite energy planes as $N \to \infty$. The
limiting finite energy planes can be chosen to foliate three
symplectic manifolds with cylindrical ends, namely the completion
$W$ of the complement of $M$ in $S^2 \times S^2$ with an end
symplectomorphic to the negative symplectization of $\Sigma$, that
is $((-\infty,0)\times \Sigma, d(e^t \alpha))$, the completion of
$U_i$, which will be a copy of $T^* S^2$, and the completion of $M
\setminus U_i$ with two ends symplectomorphic to the positive
symplectization of $\Sigma$ and the negative symplectization of
the boundary of $U_i$.

The foliations of the completion of $U_i$ were determined in
\cite{hind}. For $U_i$ and its almost-complex structure suitably
chosen, the Reeb flow on $\p U_i$ is foliated by closed orbits,
and exactly one curve in each foliation is asymptotic to each
closed orbit. Also, each curve in the foliation from ${\cal F}_0$
intersects in a single point each curve in the foliation from
${\cal F}_1$ provided that the curves have different asymptotic
limits. By the positivity of intersections, any intersections of
limiting finite energy planes must also be seen as intersections
of holomorphic spheres in the foliations ${\cal F}_0$ and ${\cal
F}_1$. It follows that finite energy curves in the foliations of
$W$ and the completion of $M\setminus U_i$ either have the same
image or are disjoint. For, if any of these curves were to
intersect in an isolated point we could find $J_{N,i}$-holomorphic
curves in ${\cal F}_0$ and ${\cal F}_1$, for $N$ sufficiently
large, with intersection number at least two. There would be one
intersection point near our point in $W$ or $M\setminus U_i$ and
another inside $U_i$. This gives a contradiction. Thus in
particular the two foliations of $W$ must coincide. Two curves in
different homotopy classes which are disjoint in the completion of
$M$ come from taking a limit of $J_{N,i}$-holomorphic curves
through the same point of $\Delta$. The resulting finite energy
curves in $W$ coincide.

We can also obtain finite energy foliations of $W$ and the
completion of $M$ by stretching the neck only along $\Sigma$ and
studying limits of $J_{N,t}$-holomorphic curves. A further limit
of the $J_{N,0}$ and $J_{N,1}$ finite energy foliations of $M$ as
we degenerate the almost-complex structure along $\p U_i$ for
$i=0,1$ would give our foliations of completions of $M \setminus
U_i$ and $U_i$ as above. If we take a diagonal subsequence of $N
\to \infty$ then for a countable dense subset of $t \in [0,1]$ the
corresponding foliations of $W$ and $M$ can be assumed to arise
from the same subsequence of $N \to \infty$. Since the $J_{N,t}$
are all equal outside of $M$, for each $t$ we obtain a finite
energy foliation of $W$ with respect to the same almost-complex
structure. For any fixed $N$, as $t \to t_0$, the
$J_{N,t}$-holomorphic foliations of $S^2 \times S^2$ converge to
the $J_{N,t_0}$-holomorphic foliations. Taking a further limit, we
see that for $t$ in our dense subset, and hence for all $t$ since
the subset was arbitrary, the finite energy foliations of $W$ vary
continuously with $t$. Suppose that the contact form $\alpha$ on
$\Sigma$ and the almost-complex structure on $W$ are chosen
generically so that periodic orbits of the Reeb flow on $\Sigma$
are isolated and embedded finite energy curves in $W$ appear in
families whose dimension is as predicted by the index theorem, see
\cite{egh} and \cite{hoff}. Then the following is true.

\begin{lemma}
The deformation index of generic finite energy curves $C$ in the foliation of $W$ satisfies $\mathrm{index} (C) \le 2$.
\end{lemma}

Since $2$-dimensional families of curves are needed to foliate
$W$, this lemma implies that all finite energy planes sufficiently
close to a curve $C$ in our foliation actually appear in the
foliation. In particular, as $t$ varies the corresponding two
finite energy foliations of $W$ remain the same. We will use this
result to draw our conclusion about the Lagrangian isotopy classes
and then prove Lemma $8$ at the end of the section.

Another result coming from the analysis in \cite{hind} is that each curve in the foliations
of the completion of $U_i$  coming from ${\cal F}_0$ and ${\cal
F}_1$ intersects $L_i$ transversally in a single point. Hence, for
$N$ sufficiently large the $J_{N,i}$ holomorphic curves in ${\cal
F}_0$ and ${\cal F}_1$ must also intersect $L_i$ transversally.

We now study the foliations ${\cal F}_0$ and ${\cal F}_1$
corresponding to the family of almost-complex structures $J_{N,t}$
on $S^2 \times S^2$ for large $N$. As in section $3$ we can find corresponding
diffeomorphisms $f_{N,t}:\Delta \to \Delta$ such that
$\overline{\Delta}=L_0$ consists of the intersections of
$J_{N,0}$-holomorphic spheres in ${\cal F}_0$ through points $x\in
\Delta$ with the spheres in ${\cal F}_1$ through $f_{N,0}(x)$, and
$L=L_1$ consists of the intersections of $J_{N,1}$-holomorphic
curves in ${\cal F}_0$ through points $x\in \Delta$ with the
spheres in ${\cal F}_1$ through $f_{N,1}(x)$. By the result in
section $2$ we may assume that $f_{N,t}(x) \neq x$ for all $x$,
$t$. As $N \to \infty$, the convergence of $J_{N,0}$ and $J_{N,1}$-holomorphic spheres implies that we may assume that the $f_{N,0}$ and $f_{N,1}$ converge to continuous maps $f_{\infty,0}$ and $f_{\infty,1}$ and hence that the $f_{N,t}$ converge to maps $f_{\infty,t}$ satisfying $f_{\infty,t}(x) \neq x$ for all $t$ and $x$.

Next we define a family of embedded spheres by letting $L_t$ be
the union of the intersections of $J_{N,t}$-holomorphic curves in
${\cal F}_0$ through points $x$ with $J_{N,t}$-holomorphic curves
in ${\cal F}_1$ through $f_{N,t}(x)$. We claim that for $N$
sufficiently large $L_t \subset M$ for all $t$. For otherwise
there exists a subsequence $N \to \infty$ such that for each $N$
there exists a $t$ and $x\in \Delta$ with the $J_{N,t}$
holomorphic sphere in ${\cal F}_0$ through $x$ intersecting the
$J_{N,t}$ holomorphic sphere in ${\cal F}_1$ through
$f_{N,t}(x)\neq x$ in a point outside of $M$. Taking a limit as $N
\to \infty$ we then have that the two foliations of $W$
corresponding to limits of $J_{N,t}$-holomorphic spheres in ${\cal
F}_0$ and ${\cal F}_1$ will not coincide, a contradiction since
all of these foliations are identical for all $t$.

Following the method of section $3$, we can find a family of
symplectic forms $\omega_t$ on $S^2 \times S^2$ such that $L_t$ is
Lagrangian with respect to $\omega_t$. The $\omega_t$ restrict to
exact symplectic forms on $M$, say $\omega_t =d \alpha_t$ which
are tamed by $J_t$. In a tubular neighbourhood
$V=(-\epsilon,0)\times \Sigma$ of the boundary $\Sigma =
\{0\}\times \Sigma$ of $M$, define a function $\chi :V \to [0,1)$
such that $\chi(r,y)$ is an increasing function of $r$,
$\chi(r,y)=0$ for $r$ close to $-\epsilon$ and $\chi(r,y)=1$ for
$r$ close to $0$. Then, first scaling $\alpha_t$ if necessary, we
can replace it by $\beta_t=(1-\chi)\alpha_t + \chi e^r \alpha$ in
$V$. The new form $\omega_t = d\beta_t$ will still be symplectic
and tamed by $J_t$ (for $\alpha_t$ suitably scaled) but now agrees
with $\omega$ near $\Sigma$. Assuming $V$ to be disjoint from all
$L_t$, the submanifolds $L_t$ will still by Lagrangian with
respect to $\omega_t$.

We now apply Moser's method as in section $3$ to find a
symplectomorphism between $(M,\omega_t)$ and $(M,\omega)$ and
thereby isotope the $L_t$ into Lagrangian submanifolds of
$(M,\omega)$. As before, this can be arranged to fix $L_0$ and
$L_1$ and now also the neighbourhood $V$. Thus it gives our
Lagrangian isotopy as required.

{\bf Proof of Lemma $8$}

Let $C$ be a finite energy curve in the foliation of $W$. The
curve $C$ will be one component of a limit of holomorphic spheres.
The other components can be assumed to curves $D_i$ in the
symplectization of $\Sigma$ and curves $E_j$ in the completion of
$M$. Since the finite energy curves are limits of curves of genus
$0$ and as in \cite{hind} the limiting curve has only one
component in $W$, the $D_i$ and $E_j$ have only one positive
asymptotic limit. We recall that the index formula for a finite
energy curve $F$ depends upon a trivialization of the contact
planes along its asymptotic limits (which are Reeb orbits in
certain contact manifolds). We can then define a Conley-Zehnder
index for each asymptotic limit and a Chern class $c_1(F)$
relative to these trivializations. Suppose that the positive
asymptotic limits have Conley-Zehnder indices $\mu_k^+$ for $1 \le
k \le m$ and the negative asymptotic limits have index $\mu_l^-$
for $1 \le l \le n$. If the asymptotic limits are nondegenerate
the formula for the deformation index of $F$ modulo
reparameterizations is
$$\mathrm{index}(F) = -(2-m-n)+2c_1(F)+ \sum_{k=1}^{m} \mu_k^+
-\sum_{l=1}^{n} \mu_l^-.$$ In our case, a global trivialization of
the contact planes in $T^1 S^2$ extends over any $1$-handles to a
trivialization of $\xi = \{\alpha =0\}$ on $\Sigma$ and we compute
our indices relative to this. Then the curves in the
symplectization of $\Sigma$ and the completion of $M$ have Chern
class $0$ and our curve $C$ has $c_1(C)=2$.

A curve $E_j$ has a single positive asymptotic limit. If this has
index $\mu_j^+$ then we obtain
$$\mathrm{index}(E_j)= -1 + \mu_j^+.$$
But for a generic choice of almost-complex structure this index
must be nonnegative and so $\mu_j^+ \ge 1$ for all $j$.

Suppose that a curve $D_i$ has a positive asymptotic limit with
index $\mu_i^+$ and $m_i$ negative asymptotic limits with index
$\mu_{ik}^-$ for $1 \le k \le m_i$. Then the index formula becomes
$$\mathrm{index}(D_i)= -1 +m_i + \mu_i^+
-\sum_{k=1}^{m_i}\mu_{ik}^-.$$ Again generically this index must
be nonnegative. Therefore if all of the negative asymptotic limits
are positive asymptotic limits of curves $E_j$ we obtain $\mu_i^+
\ge 1$. By an induction on the number of levels of the limiting
finite energy curve we deduce that in fact $\mu_i^+ \ge 1$ for all
curves $D_i$.

Finally we look at $C$. It has only negative asymptotic limits and
if these have index $\mu_l^-$ for $1\le l \le n$ then
$$\mathrm{index}(C) = 2+n-\sum_{l=1}^{n}\mu_l^-.$$
But all of these negative asymptotic limits are positive limits of
curves $D_i$ or $E_j$. Hence $\mu_l^- \ge 1$ for all $l$ and so
$\mathrm{index}(C) \le 2$ as required.

\end{document}